\documentclass[11pt]{amsart}
\usepackage{latexsym,amssymb,amsmath,amsfonts,bm}
\usepackage{amssymb}
\usepackage{amsmath}
\usepackage{graphicx}
\usepackage{dcolumn}
\usepackage{hyperref}
\usepackage[utf8]{inputenc}
\usepackage[english]{babel}
\usepackage{amsthm}
\usepackage{latexsym}
\usepackage{xcolor}
\usepackage{flexisym}
\usepackage{ulem}

\usepackage{subfig}
\usepackage{float}

\newtheorem{theorem}{Theorem}[section]
\newtheorem{corollary}[theorem]{Corollary}

\newtheorem{remark}[theorem]{Remark}
\newtheorem{proposition}[theorem]{Proposition}
\newtheorem{definition}[theorem]{Definition}
\newtheorem{example}{Example}[section]

\newtheorem{remarks}[theorem]{Remarks}

\newtheorem{conjecture}[theorem]{Conjecture}

\newcommand{\boldalpha}{\boldsymbol\alpha}

\newlength{\cellsize}
\newcommand\tableau[1]{
\vcenter{
\let\\=\cr
\baselineskip=-16000pt
\lineskiplimit=16000pt
\lineskip=0pt
\halign{&\tableaucell{##}\cr#1\crcr}}}

\newcommand{\tableaucell}[1]{{%
\def \arg{#1}\def \void{}%
\ifx \void \arg
\vbox to \cellsize{\vfil \hrule width \cellsize height 0pt}%
\else
\unitlength=\cellsize
\begin{picture}(1,1)
\put(0,0){\makebox(1,1){$#1$}}
\put(0,0){\line(1,0){1}}
\put(0,1){\line(1,0){1}}
\put(0,0){\line(0,1){1}}
\put(1,0){\line(0,1){1}}
\end{picture}%
\fi}}



\usepackage{lipsum}

\usepackage[left=1in, right=1in, top=1in, bottom=1in]{geometry}

\begin{document}

\title[Towards a Combinatorial Model for $q$-weight Multiplicities of Simple Lie Algebras]{Towards a Combinatorial Model for $q$-weight Multiplicities of Simple Lie Algebras (Extended Abstract)}

\author[C. Lecouvey]{C\'{e}dric Lecouvey}
\address[C\'{e}dric Lecouvey]{Laboratoire de Math\'{e}matiques et Physique Th\'{e}orique, Universit\'{e} de Tours, Tours, Centre-Val de Loire, 37020, France}
\email{cedric.lecouvey@lmpt.univ-tours.fr}

\author[C.~Lenart]{Cristian Lenart}
\address[Cristian Lenart]{Department of Mathematics and Statistics, State University of New York at Albany, 
Albany, NY 12222, U.S.A.}
\email{clenart@albany.edu}

\author[A.~Schultze]{Adam Schultze}
\address[Adam Schultze]{Department of Mathematics, Statistics, and Computer Science, St. Olaf College, 1520 St. Olaf Avenue, 
Northfield, MN 55057, U.S.A.}
\email{schult24@stolaf.edu}


%% you can include your bibliography however you want, but using an external .bib file is STRONGLY RECOMMENDED and will make the editor's life much easier
%% regardless of how you do it, please use numerical citations, ie. [xx, yy] in the text

%% this sample uses biblatex, which (among other things) takes care of URLs in a more flexible way than bibtex
%% but you can use bibtex if you want
%\usepackage[backend=bibtex]{biblatex}  This line
%\addbibresource{sample.bib}            and this line (undo for original setup)
%% note the \printbibliography command at the end of the file which goes with these biblatex commands

%% note that you DO NOT have to put your abstract here -- it is generated by \maketitle and the \abstract and \resume commands above

\begin{abstract}
Kostka-Foulkes polynomials are Lusztig's $q$-analogues of weight multiplicities for irreducible representations of semisimple Lie algebras. It has long been known that these polynomials have non-negative coefficients.  A statistic on semistandard Young tableaux with partition content, called \textit{charge}, was used to give a combinatorial formula exhibiting this fact in type $A$. Defining a charge statistic beyond type $A$ has been a long-standing problem. Here, we take a completely new approach based on the definition of Kostka-Foulkes polynomials as an alternating sum over Kostant partitions, which can be thought of as formal sums of positive roots. We use a sign-reversing involution to obtain a positive expansion, in which the relevant statistic is simply the number of parts in the Kostant partitions.  The hope is that the simplicity of this new crystal-like model will naturally extend to other classical types.
\end{abstract}

\maketitle

%\vspace{12pt}\textit{Abstract:} 

%\newpage
%\tableofcontents
%\newpage
\section{Introduction}
 It has long been known that the Kostka-Foulkes polynomials (that is, Lusztig's $q$-analogues of weight multiplicities) have non-negative coefficients \cite{Lusztig}.  A statistic on semistandard Young tableaux (SSYT) with partition content, called \textit{charge}, was used to give a combinatorial formula exhibiting this fact in type $A$ \cite{LS}. Defining a charge statistic beyond type $A$ has been a long-standing problem. However, Lecouvey was able to define such a statistic on Kashiwara-Nakashima tableaux of classical Lie type \cite{Lecouvey 2005,Lecouvey 2006}, but was only able to relate his charge to the corresponding Kostka-Foulkes polynomials in special cases.

 In \cite{Lascoux}, Lascoux gave another combinatorial formula to show the non-negativity of Kostka-Foulkes polynomials in type $A$ by decomposing them into so-called atomic polynomials. Unfortunately, the complexity of the decomposition made it difficult to extend to other Lie types. Recently, Lenart and Lecouvey gave a new model for such an atomic decomposition with a much simpler combinatorial structure, allowing them to extend beyond type $A$.  In \cite{Lec Lenart}, they explicitly described this decomposition in all classical types for the stable cases (meaning for sufficiently large rank).  However, it has proven difficult to relate the atomic decomposition to the charge statistic beyond type $A$.

 Here, we take a completely different approach based on the definition of Kostka-Foulkes polynomials as an alternating sum over Kostant partitions as opposed to the above models which are based on tableaux. We combine the structure of the crystals $B(\infty)$ and $B^*(\infty)$ to construct a collection of modified crystal operators on Kostant partitions.  These modifications are indeed needed, as the crystal structures $B(\infty)$ and $B^*(\infty)$ on their own provide too few arrows which preserve the partition size, so we cannot use them individually to cancel all the negative terms in the alternating formula. On the other hand, we note that the classical embedding of SSYT into $B(\infty)$ does not give way to the Kostka-Foulkes polynomials as the generating function of the number of parts of Kostant partitions. We then define new operators on this modified crystal graph which act similarly to the Weyl action on classical crystals.  The resulting graphs have some remarkable properties: they satisfy the braid relations (although with occurrences of vertices with value zero) and they are subgraphs of the Cayley graphs of parabolic subgroups of $S_n$.  We then provide a matching on these graphs (based on two computer verified conjectures, both of which pertain to certain occurences of braid relations which contain vertices of value zero) which corresponds to a sign-reversing involution that cancels out all of the negative terms in the corresponding Kostka-Foulkes polynomial.  This positive expansion in terms of Kostant partitions gives way to a statistic which is simply read by counting the number of parts of the Kostant partitions.   The hope is that the simplicity of this new crystal-like model will naturally extend to other classical types.

%\section{Background}
%\subsection{Kostka-Foulkes Polynomials and the Kostant Partition Function}
\section{Kostka-Foulkes Polynomials and the Kostant Partition Function}
For a finite simple Lie algebra $\mathfrak{g}$, let $\Phi$ be a finite root system with simple roots $\alpha_i$ and $W$ be its corresponding Weyl group with $\rho$ being the half sum of all positive roots.  Give $\Phi^+$, the set of positive roots, a total ordering so that we may enumerate the positive roots as $\boldsymbol\alpha_i$ for $i\in \{1,2,\hdots,N\}$ where $N = |\Phi^+|$.  For dominant weight $\lambda$ and a weight $\mu$, let $K_{\lambda,\mu}$, commonly referred to as \textit{Kostka numbers}, be the multiplicity of the weight $\mu$ in the irreducible finite dimensional representation $V(\lambda)$ of $\mathfrak{g}$.  The next definitions follow closely to the notation given in \cite{Lee Salisbury}.

%\vspace{12pt}
  Set $K = \{(\beta):\beta\in\Phi^+\}$, where each $(\beta)$ is considered as a formal symbol.  Let $K\textprime{}$ be the free abelian group generated by $K$ and then define $\mathcal{K}$ to be the set of elements in $K\textprime{}$ with coefficients from $\mathbb{N}$.  An element of $\mathcal{K}$ will then be considered as a Kostant partition, and will be denoted by a boldface Greek letter.  For $\boldalpha\in\mathcal{K}$, we have $\boldalpha = \displaystyle\sum\limits^N_{i=1}c_i(\beta_i)$ and we then define the \textit{evaluation} of $\boldalpha$ to be $\overline{\boldalpha} = \displaystyle\sum\limits^N_{i=1}c_i\beta_i$.

%\vspace{12pt} 
The Kostka numbers can be written in terms of the Kostant partition function, $\mathcal{P},$ which is given by the equality
$\displaystyle\prod_{\beta\in \Phi^+} \dfrac{1}{1-e^{\beta}} = \displaystyle\sum_{\gamma\in\mathcal{K}}\mathcal{P}(\gamma)e^{\gamma}.$  It is worth noting that $\mathcal{P}$ then counts the number of ways the weight $\gamma$ can be written as a sum of positive roots.  The Weyl character formula then gives us the identity $$K_{\lambda,\mu} = \displaystyle\sum_{w\in W} sgn(w)\mathcal{P}(w(\lambda + \rho) - (\mu + \rho)).$$

%\begin{equation}\label{kostka number in terms of kostant partition function}
%K_{\lambda,\mu} = \sum_{w\in W} sgn(w)\mathcal{P}(w(\lambda + \rho) - (\mu + \rho)).
%\end{equation}

We then generalize this to define the polynomial $K_{\lambda,\mu}(t)$. These polynomials are often referred to as Lusztig's \textit{$t$-analogue} of weight multiplicity.  We do so by defining a $t$-analogue for the Kostant partition function, $\mathcal{P}_t$, given by the equality $$\prod_{\beta\in \Phi^+} \dfrac{1}{1-te^{\beta}} = \sum_{\gamma\in\mathcal{K}}\mathcal{P}_t(\gamma)e^{\gamma}.$$  Here, $\mathcal{P}_t$ is a polynomial in $t$ where, for $\gamma$ in the associated root lattice, the coefficient of $t^l$ is the number of ways to partition $\gamma$ into an unordered sum of $l$ positive roots.  It is then useful to explicitly write
$\mathcal{P}_t(\gamma) := \displaystyle\sum_{\substack{(n_1,\hdots,n_N)\in\mathbb{N}^N \\ n_1\beta_1 + \hdots + n_N\beta_N = \gamma}} t^{n_1+\hdots + n_N}.$

%\noindent and then $K_{\lambda,\mu}(t) := \sum_{w\in W} sgn(w)\mathcal{P}_t(w(\lambda + \rho) - (\mu + \rho)).$

%\begin{equation}
%\mathcal{P}_t(\gamma) := \sum_{\substack{(n_1,\hdots,n_N)\in\mathbb{N}^N \\ n_1\beta_1 + \hdots + n_N\beta_N = \gamma}} t^{n_1+\hdots + n_N}.
%\end{equation}

%\begin{definition}
%The Kostka-Foulkes polynomials $K_{\lambda,\mu}(t)$ are given by
%$$K_{\lambda,\mu}(t) := \sum_{w\in W} sgn(w)\mathcal{P}_t(w(\lambda + \rho) - (\mu + \rho)).$$

%\begin{equation}
%K_{\lambda,\mu}(t) := \sum_{w\in W} sgn(w)\mathcal{P}_t(w(\lambda + \rho) - (\mu + \rho)).
%\end{equation}
%\end{definition}

\begin{definition}
{\rm
We define Lusztig's $t$-analogue of weight multiplicity as
$$K_{\lambda,\mu}(t) := \displaystyle\sum_{w\in W} sgn(w)\mathcal{P}_t(w(\lambda + \rho) - (\mu + \rho)).$$

%\begin{equation}
%K_{\lambda,\mu}(t) := \sum_{w\in W} sgn(w)\mathcal{P}_t(w(\lambda + \rho) - (\mu + \rho)).
%\end{equation}
}
\end{definition}

When $\mu$ is also a dominant weight and $\lambda - \mu$ is the positive sum of simple roots, the polynomials $K_{\lambda,\mu}(t)$ have positive coefficients and are known as Kostka-Foulkes polynomials.

The equation defining $K_{\lambda,\mu}(t)$ can be greatly simplified using the notation of Kostant partitions.  Let $\boldsymbol\alpha = \sum\limits_{i=1}^N c_i(\beta_i)$ be a Kostant partition.  Then the \textit{size} of $\boldsymbol\alpha$ is given by $|\boldsymbol\alpha| =\sum\limits_{i=1}^N c_i$.  This allows us to rewrite the $t$-analogue Kostant partition function as 
$\mathcal{P}_t(\gamma) = \displaystyle\sum\limits_{\boldsymbol\alpha\in\mathcal{K}, \overline{\boldsymbol\alpha} = \gamma} t^{|\boldsymbol\alpha|}.$

%\begin{equation}
%\mathcal{P}_t(\gamma) = \sum\limits_{\boldsymbol\alpha\in\mathcal{K}, \overline{\boldsymbol\alpha} = \gamma} t^{|\boldsymbol\alpha|}.
%\end{equation}

If we then let 
$S_{\lambda,\mu} := \{(w,\boldsymbol\alpha): w\in W, \boldsymbol\alpha\in\mathcal{K}\hspace{4pt}\text{where} \hspace{4pt} \overline{\boldsymbol\alpha} = w(\lambda + \rho) - (\mu + \rho)\},$
%\begin{equation}
%S_{\lambda,\mu} := \{(w,\boldsymbol\alpha): w\in W, \boldsymbol\alpha\in\mathcal{K}\hspace{4pt}\text{where} \hspace{4pt} \overline{\boldsymbol\alpha} = w(\lambda + \rho) - (\mu + \rho)\},
%\end{equation}
we can rewrite the Kostka-Foulkes polynomials as
\begin{equation}
K_{\lambda,\mu}(t) = \sum_{w\in W}\left( sgn(w) \sum_{(w,\boldsymbol\alpha)\in S_{\lambda,\mu}} t^{|\boldsymbol\alpha|}\right) = \sum_{(w,\boldsymbol\alpha)\in S_{\lambda,\mu}}sgn(w) t^{|\boldalpha|}.
\end{equation}

\vspace{12pt}
 This paper will focus on type $A_{n-1}$ so let us start with the basic facts of this type.  We can identify the space $\mathfrak{h}^*_{\mathbb{R}}$ with the quotient $V:=\mathbb{R}^n/\mathbb{R}(1,\ldots,1)$, where $\mathbb{R}(1,\ldots,1)$ denotes the subspace in $\mathbb{R}^n$ spanned by the vector $(1,\ldots,1)$.  Let $\varepsilon_1,\ldots,\varepsilon_n\in V$ be the images of the coordinate vectors in $\mathbb{R}^n$.  The root system is $\Phi = \{\alpha_{ij} := \varepsilon_i-\varepsilon_j : i\neq j, 1\leq i,j \leq n\}$. The simple roots are $\alpha_i = \alpha_{i,i+1}$, for $i = 1,\ldots,n-1$ and the positive roots are $\Phi^+=\{\alpha_{ij} := \varepsilon_i-\varepsilon_j : i\neq j, 1\leq i<j \leq n\}$. It will be beneficial for us to think of the positive roots as tuples with the notation $(i,j)$ for the root $\alpha_{i,j}$.  The weight lattice is $\Lambda = \mathbb{Z}^n/\mathbb{Z}(1,\ldots,1)$.  The fundamental weights are $\omega_i = \varepsilon_1 + \varepsilon_2 + \ldots + \varepsilon_i$, for $i = 1,2,\ldots,n-1$. A dominant weight $\lambda = \lambda_1\varepsilon_1 + \ldots + \lambda_{n-1}\varepsilon_{n-1}$ is identified with the partition $(\lambda_1\geq\lambda_2\geq\ldots\geq\lambda_{n-1}\geq\lambda_n=0)$ having at most $n-1$ parts.  Note that $\rho = (n-1,n-2,\ldots,0)$.  We now produce an example of $S_{\lambda,\mu}$ and $K_{\lambda,\mu}(t)$ in type $A_{n-1}$.

\begin{example}\label{example KF poly}
{\rm We now build the Kostka Foulkes polynomial $K_{\lambda,\mu}(t)$ of type $A_3$ for \\$\lambda = (2,2,0,0)$ and $\mu = (1,1,1,1)$. First, we build the associated $S_{\lambda,\mu}$.  Recall that \\$S_{\lambda,\mu} = \{(w,\boldsymbol\alpha): w\in W, \boldsymbol\alpha\in\mathcal{K}\hspace{4pt}\text{where} \hspace{4pt} \overline{\boldsymbol\alpha} = w(\lambda + \rho) - (\mu + \rho)\}.$  We then get the following:

%We have the following weights for $w\in W$:
%$$(\lambda + \rho) - (\mu + \rho) = (1,1,-1,-1)$$
%$$s_1(\lambda + \rho) - (\mu + \rho) = (0,2,-1,-1)$$
%$$s_2(\lambda + \rho) - (\mu + \rho) = (1,-2,2,-1)$$
%$$s_3(\lambda + \rho) - (\mu + \rho) = (1,1,-2,0)$$
%$$s_2s_1(\lambda + \rho) - (\mu + \rho) = (0,-2,3,-1)$$
%$$s_1s_2(\lambda + \rho) - (\mu + \rho) = (-3,2,2,-1)$$
%$$s_3s_1(\lambda + \rho) - (\mu + \rho) = (0,2,-2,0)$$
%$$s_1s_2s_1(\lambda + \rho) - (\mu + \rho) = (-3,1,3,-1)$$
%$$\vdots$$
%Note that only $1$, $s_1$, $s_3$, and $s_3s_1 = s_1s_3$ are the only words which give way to a weight of a Kostant partition (the others each require a negative coefficient of a positive root to obtain the desired weight).  We then get the following:
\[S_{(2,2,0,0),(1,1,1,1)} = \{[1,(\alpha_1 + \alpha_2) + (\alpha_2+\alpha_3)],[1,(\alpha_1) + (\alpha_2) + (\alpha_2+\alpha_3)], [1, (\alpha_1 + \alpha_2+\alpha_3) + (\alpha_2)],\] \[[1, (\alpha_1+\alpha_2) + (\alpha_2) + (\alpha_3)], [1, (\alpha_1) + 2(\alpha_2) + (\alpha_3)], [s_1, (\alpha_2) + (\alpha_2+\alpha_3)],\] \[ [s_1, 2(\alpha_2) + (\alpha_3)], [s_3,(\alpha_1+\alpha_2)+(\alpha_2)], [s_3, (\alpha_1) + 2(\alpha_2)],\] \[ [s_3s_1, 2(\alpha_2)] \}.\] 
Then $$K_{(2,2,0,0),(1,1,1,1)}(t) = \sum_{(w,\boldsymbol\alpha)\in S_{(2,2,0,0),(1,1,1,1)}}sgn(w) t^{|\boldalpha|} $$
\vspace{8pt}
 $$= t^2 + t^3 + t^2 + t^3 + t^4 - t^2 - t^3 - t^2 - t^3 + t^2 = t^4 + t^2.$$

\vspace{12pt}In this paper we will provide a combinatorial method, through the use of crystals and Kostant partitions, to match like terms of opposite signs.  The unmatched terms will then give way to the desired positive form of the Kostka-Foulkes polynomials.
}
\end{example}

%\subsection{The Crystals $B(\infty)$ and $B^*(\infty)$}
\section{The Crystals $B(\infty)$ and $B^*(\infty)$}\label{Section infty crystals}

We can endow the set $\mathcal{K}$ of type $A_{n-1}$ with two different crystal structures: $B(\infty)$, which comes from $U_q^-(\mathfrak{sl}_{n})$, and $B^*(\infty)$, the image of $B(\infty)$ by the Kashiwara $*$-involution. For more detail on these crystals refer to \cite{Bump and Schilling}.  This section describes realizations of $B(\infty)$ and $B^*(\infty)$ on Kostant partitions as developed in \cite{claxting}. 
We first define two orderings on the positive roots which will be used for computing the crystal operators on a given Kostant partition $\boldalpha\in \mathcal{K}$.  

\vspace{12pt}
\begin{definition}\label{words of string definitions}
{\rm
Given a Kostant partition $\boldalpha\in \mathcal{K}$, define $W_i^{\boldalpha}$ to be the string as follows:
\begin{enumerate}
\item Order the roots $(a,b)$ of $\boldalpha$ from left to right first decreasing by $a$, and for the same $a$, increasing by $b$.

%\item Order the roots $(h,j)$ of $\boldalpha$ from left to right, first in increasing order of height, then by largest to smallest $h$ value.
\item Place a  $+$ below each root $(h, i)$ and a $-$ below each root $(h, i+1)$.
\end{enumerate}
Then define $W_i^{*,\boldalpha}$ to be the string with the following rule:
\begin{enumerate}
\item Order the roots $(a,b)$ of $\boldalpha$ from left to right first increasing by $b$, and for the same $b$, decreasing by $a$.
%\item Order the roots $(h,j)$ of $\boldalpha$ from left to right, first in increasing order of height, then by smallest to largest $h$ value.
\item Place a  $+$ below each root $(i+1, j)$ and a $-$ below each root $(i, j)$.
\end{enumerate}

\vspace{8pt}
We then define the reduced words $W_i^{\boldalpha,r}$ and $W_i^{*,\boldalpha,r}$ to be the result of canceling out any pairs of $-$ and $+$ where the $-$ appears directly to the left of  a $+$ or if the only symbols between them have been previously canceled in the words $W_i^{\boldalpha}$ and $W_i^{*,\boldalpha}$ respectively. This so-called \textit{$i$-pairing procedure} will be used again in Section~\ref{crystal structure on KP} when we define a similar word. For the remainder of the paper we will drop the $\boldalpha$ superscript unless there is ambiguity between Kostant partitions. }
\end{definition}

%\begin{remark}\label{ipair}{\rm Given a $\pm$ word $W_i$, let $W_i^r=(+)^a(-)^b$ be the result of the pairing procedure mentioned above. An easy way to visualise this procedure is the following: consider the piecewise-linear graph starting at $(0,0)$, which is obtained from $W_i$ by associating an up step $ (x,y) \rightarrow(x+1,y+1)$ to each $+$, and a down step $ (x,y) \rightarrow(x+1,y-1)$ to each $-$.  By viewing the graph as a sequence of mountains and valleys (obscuring each other), the letters $+$ in $W_i^r$ correspond to the up steps visible from the left, while the letters $-$ in $W_i^r$ correspond to the down steps visible from the right. In particular, $a$ is the maximum height of the graph, while $b$ is the difference between the maximum height and the height of the endpoint of the graph. Moreover, the rightmost $+$ (respectively leftmost $-$) in $W_i^r$ corresponds to the step before the first maximum (respectively after the last maximum) of the graph.
%}
%\end{remark}

\begin{definition}\label{Binfty operators}
{\rm For $\boldalpha\in\mathcal{K}$ and $1\leq i\leq n-1$, define

 \[\widetilde{f}_i(\boldalpha) := \begin{cases} 
      \boldalpha - (h,i) + (h,i+1)  & \text{if the rightmost + in} \hspace{4 pt} W_i^{r}\hspace{4 pt} \text{is below an} \hspace{4pt} (h,i) \\
      \boldalpha + (i,i+1)  & \text{if there are no + in }\hspace{4pt} W_i^{r}
   \end{cases}\]
   
   \[\widetilde{e}_i(\boldalpha) := \begin{cases} 
      \boldalpha - (h,i+1) + (h,i)  & \text{if the leftmost $-$ in} \hspace{4 pt} W_i^{r}\hspace{4 pt} \text{is below an} \hspace{4pt} (h,i+1), h\neq i \\
      \boldalpha - (i,i+1)  & \text{if the leftmost $-$ in} \hspace{4 pt} W_i^{r}\hspace{4 pt} \text{is below an} \hspace{4pt} (i,i+1)  \\
      0  & \text{if there are no $-$ in }\hspace{4pt} W_i^{r}
   \end{cases}\]

    \[\widetilde{f}^*_i(\boldalpha) := \begin{cases} 
      \boldalpha - (i+1,j) + (i,j)  & \text{if the rightmost + in} \hspace{4 pt} W_i^{*,r}\hspace{4 pt} \text{is below an} \hspace{4pt} (i+1,j) \\
      \boldalpha + (i,i+1)  & \text{if there are no + in }\hspace{4pt} W_i^{*,r}
   \end{cases}\]
   
   \[\widetilde{e}^*_i(\boldalpha) := \begin{cases} 
      \boldalpha - (i,j) + (i+1,j)  & \text{if the leftmost $-$ in} \hspace{4 pt} W_i^{*,r}\hspace{4 pt} \text{is below an} \hspace{4pt} (i,j), j\neq i+1 \\
      \boldalpha - (i,i+1)  & \text{if the leftmost $-$ in} \hspace{4 pt} W_i^{*,r}\hspace{4 pt} \text{is below an} \hspace{4pt} (i,i+1)  \\
      0  & \text{if there are no $-$ in }\hspace{4pt} W_i^{*,r}
   \end{cases}\]
   
       \[\widetilde{\varphi}_i(\boldalpha) := a_i \hspace{36pt} \widetilde{\varepsilon}_i(\boldalpha) := b_i \hspace{36pt} \widetilde{\varphi}^*_i(\boldalpha) = a^*_i \hspace{36pt} \widetilde{\varepsilon}^*_i(\boldalpha) = b^*_i.\]
}
\end{definition}

\begin{proposition}\cite{claxting}\label{bicrystal and wt prop}
The set $\mathcal{K}$ along with the crystal operators $\widetilde{f}_i$, $\widetilde{e}_i$ and string operators $\widetilde{\varphi}_i$, $\widetilde{\varepsilon}_i$  give way to the crystal $B(\infty)$ and the set $\mathcal{K}$ along with the crystal operators $\widetilde{f}_i^*$, $\widetilde{e}^*_i$ and string operators $\widetilde{\varphi}^*_i$, $\widetilde{\varepsilon}^*_i$ give way to the crystal $B^*(\infty)$. 

\noindent Further, we have that:
$\langle wt(\boldalpha), \alpha_i^{\vee}\rangle = a_i + a_i^{*} - b_i - b_i^{*}.$\end{proposition}
%\vspace{-12pt}\begin{equation}\label{wt equation def}
%\langle wt(\boldalpha), \alpha_i^{\vee}\rangle = a_i + a_i^{*} - b_i - b_i^{*}.
%\end{equation}
%\end{proposition}

 These are \textit{Stembridge crystals}, meaning that they adhear to the following axioms \cite{Bump and Schilling}.  Each of these axioms have dual forms ($S0\textprime{}-S3\textprime{}$) in terms of the $f_i$ operators. 

\vspace{12pt}\noindent\textbf{Stembridge Axiom S0.} If $e_i(x) = 0$, then $\varepsilon_i(x) = 0$.

%\vspace{12pt}\textbf{Stembridge Axiom S0\textprime{}.} If $f_i(x) = 0$, then $\varphi_i(x) = 0$.

\vspace{12pt}
\noindent\textbf{Stembridge Axiom S1.} When $i,j\in I$ and $i\neq j$, if $x,y\in B$ and $y = e_ix$, then $\varepsilon_j(y)$ equals either $\varepsilon_j(x)$ or $\varepsilon_j(x)+1$.  The second case where $\varepsilon_j(y) = \varepsilon_j(x)+1$ is only possible when $\alpha_i$ and $\alpha_j$ are not orthogonal roots.

%\vspace{12pt}\textbf{Stembridge Axiom S1\textprime{}.} When $i,j\in I$ and $i\neq j$, if $x,y\in B$ and $y = f_ix$, then $\varphi_j(y)$ equals either $\varphi_j(x)$ or $\varphi_j(x)+1$.  The second case where $\varphi_j(y) = \varphi_j(x)+1$ is only possible when $\alpha_i$ and $\alpha_j$ are not orthogonal roots.

\vspace{12pt}
\noindent\textbf{Stembridge Axiom S2.} Assume that $i,j\in I$ and $i\neq j$.  If $x\in B$ with $\varepsilon_i(x)>0$ and $\varepsilon_j(e_ix) = \varepsilon_j(x)>0$, then $e_ie_jx = e_je_ix$ and $\varphi_i(e_jx) = \varphi_i(x)$.

%\vspace{12pt}\textbf{Stembridge Axiom S2\textprime{}.} Assume that $i,j\in I$ and $i\neq j$.  If $x\in B$ with $\varphi_i(x)>0$ and $\varphi_j(f_ix) = \varphi_j(x)>0$, then $f_if_jx = f_jf_ix$ and $\varepsilon_i(f_jx) = \varepsilon_i(x)$.

\vspace{12pt}
\noindent\textbf{Stembridge Axiom S3.} Assume that $i,j\in I$ and $i\neq j$.  If $x\in B$ with $\varepsilon_j(e_ix) = \varepsilon_j(x) + 1>1$ and $\varepsilon_i(e_jx) = \varepsilon_i(x)+ 1>1$, then 
$e_je_i^2e_jx = e_ie_j^2e_ix \neq 0,$
$\varphi_i(e_jx) = \varphi_i(e^2_je_ix)$ and $\varphi_j(e_ix) = \varphi_j(e_i^2e_jx)$.
    
%\vspace{12pt}\textbf{Stembridge Axiom S3\textprime{}.} Assume that $i,j\in I$ and $i\neq j$.  If $x\in B$ with $\varphi_j(f_ix) = \varphi_j(x) + 1>1$ and $\varphi_i(f_jx) = \varphi_i(x)+ 1>1$, then 
%\[f_jf_i^2f_jx = f_if_j^2f_ix \neq 0,\]
%$\varepsilon_i(f_jx) = \varepsilon_i(f^2_jf_ix)$ and $\varepsilon_j(f_ix) = \varepsilon_j(f_i^2f_jx)$.

%\newpage
%\section{The Main Result}\label{Section Main Result}

%\subsection{A Crystal Structure on Kostant Partitions}
\section{A Crystal Structure on Kostant Partitions}\label{crystal structure on KP}

Recall that for two dominant weights $\lambda,\mu$ where $\lambda -\mu$ is the positive sum of simple roots, we write the Kostka-Foulkes polynomials in terms of
 $ S_{\lambda,\mu} := \{(w,\boldalpha) | \overline{\boldalpha} = w (\lambda +\rho ) - (\mu + \rho )\}$
 where $\boldalpha$ is viewed as a multiset of positive roots or a Kostant partition.  In this section we will define operators that act on the Kostant partitions in a way that preserves size.  We will then use these operators to define crystal-like operators on $(w,\boldalpha)\in S_{\lambda,\mu}$ which are size preserving for $\boldalpha$ and only change the element $w\in S_n$ by the multiplication of a simple reflection on the left.  This will then be used to produce a sign-reversing size preserving involution on the elements of $S_{\lambda,\mu}$ which will give way to the main result.
 %
 %We now endow $\mathcal{K}$, the set of Kostant partitions, with the structure of a graph using the following setup.
 % Here we will drop the superscript of $\boldalpha$ on $W$ unless it is needed to differentiate between multiple words.

\vspace{12pt}Let $\boldalpha\in \mathcal{K}$ and recall the associated words $W_i^{*,r}$, $W_i^r$ from Section~\ref{Section infty crystals}. Now consider the concatenation $w_i  := W_i^{*,r}W_i^r= (+)^{a^*_i}(-)^{b^*_i}(+)^{a_i}(-)^{b_i}$ and its reduced word 
 \begin{equation}\label{reduced word def}
 w_i^r = (+)^{u_i}(-)^{v_i} = \begin{cases} 
      (+)^{a_i^*+a_i-b_i^*}(-)^{b_i} & \text{if}\hspace{5pt} b_i^* < a_i \\
      (+)^{a^*_i}(-)^{b^*_i-a_i+b_i} & \text{if}\hspace{5pt} b_i^*\geq a_i
   \end{cases}
   \end{equation}
   
We use $w_i^{r}$ to define operators that act on $\boldalpha\in \mathcal{K}$ which act similar to those given in Definition~\ref{Binfty operators}.  However, we define them in such a way that they will never change the size of the Kostant partition that they act upon.

\begin{definition}\label{modified operators def}
{\rm
For any Kostant partition $\boldalpha\in\mathcal{K}$ along with its associated word $w_i$ and $1\leq i\leq n-1$, the crystal operators $f_i,e_i$ along with the string operators $\varepsilon_i,\varphi_i$  are defined as follows:

   \[f_i(\boldalpha) := \begin{cases}
     \widetilde{f}_i(\boldalpha) & \text{if}\hspace{5pt} b_i^* <  a_i  \\
     \widetilde{f}_i^*(\boldalpha) & \text{if} \hspace{5pt} b_i^*\geq  a_i \hspace{3pt} \text{and} \hspace{5pt} a^*_i>0  \\
     0                     & \text{otherwise}
  \end{cases}\]
  
   \[e_i(\boldalpha) := \begin{cases}
     \widetilde{e}_i(\boldalpha) & \text{if}\hspace{5pt} b_i^* \le  a_i\text{, $b_i>0$, and the leftmost $-$ in} \hspace{4 pt} W_i^{r}\hspace{4 pt} \text{is not below the root} \hspace{4pt} (i,i+1) \\
     \widetilde{e}^*_i(\boldalpha) & \text{if} \hspace{5pt} b_i^*>  a_i \hspace{3pt} \text{and the leftmost $-$ in} \hspace{4 pt} W_i^{*,r}\hspace{4 pt} \text{is not below the root} \hspace{4pt} (i,i+1)  \\
     0                     & \text{otherwise}
  \end{cases}\,\]
  
  \[\widehat{e}_i(\boldalpha) := \begin{cases}
     \widetilde{e}_i(\boldalpha) & \text{if}\hspace{5pt} b_i^* \le  a_i\text{ and $b_i>0$} \\
     \widetilde{e}^*_i(\boldalpha) & \text{if} \hspace{5pt} b_i^*>  a_i  \\
     0                     & \text{otherwise}
  \end{cases}\,\]
  
  \[\varepsilon_i(\boldalpha):=\max\,\{k : e_i^k(\boldalpha)\neq 0\}, \hspace{24pt} \widehat{\varepsilon}_i(\boldalpha) := \max\,\{k : \widehat{e}_i^k(\boldalpha) \neq 0\}, \hspace{24pt} \varphi_i(\boldalpha) := \max\,\{k : f_i^k(\boldalpha)\neq 0\}\]
 }
\end{definition}

\begin{remarks}\label{rems2ef} {\rm 

(1) The definition of $f_i$ and $e_i$ on the set ${\mathcal{K}}$ of Kostant partitions can be phrased as the following composition (up to removing the arrows which add/remove roots): diagonally embed ${\mathcal{K}}$ into $B^*(\infty)\otimes B(\infty)$, then apply these operators to the tensor product (using the well-known tensor product rule), and then project onto the component in which the change occurs.

\vspace{8pt}
(2) The definitions of $f_i$, $e_i$, and  $\widehat{e}_i$ are equivalent (up to the condition on the leftmost $-$ for $e_i$) to the $i$-pairing procedure on the concatenation $W_i^{*} W_i$.

\vspace{8pt}
(3) The first two cases in the definitions of $f_i(\boldalpha)$ and $e_i(\boldalpha)$ produce a nonzero result, which is a Kostant partition with the same number of parts as $\boldalpha$. Similarly, the first two cases in the definition of $\widehat{e}_i(\boldalpha)$ also produce a nonzero result, but this Kostant partition can have one less part than $\boldalpha$, namely $\widehat{e}_i(\boldalpha)=\boldalpha-(\alpha_i)$.

   \vspace{8pt}
 (4)  Note that $\varphi_i(\boldalpha) = u_i$ but $\widehat{\varepsilon}_i, \varepsilon_i(\boldalpha) \leq v_i$ unlike classical crystals (cf. Definition~\ref{Binfty operators}).
 %.  This means that our $\varepsilon_i$ and $\widehat{\varepsilon}_i$ are not defined as they are in classical crystals.% The related implications will be addressed in Section~\ref{section stembridge proof}. 
 
 \vspace{8pt}
 (5) The weight function $wt(\boldalpha)$ is the negative sum of positive roots and should not be confused with $\overline{\boldalpha}$.  We in fact have $wt(\boldalpha) = -\overline{\boldalpha}$.
 
 \vspace{8pt}
 (6) Proposition~\ref{bicrystal and wt prop} along with Equation~\ref{reduced word def} give that $wt_i(\boldalpha):= \langle wt(\boldalpha), \alpha_i^{\vee}\rangle = u_i - v_i$.
 }
 \end{remarks}
 
 The theorems below show that the $f_{\cdot}$ and $e_{\cdot}$ operators follow the Stembridge axioms (cf. Section~\ref{Section infty crystals}) with the below possible exceptions found when applying $e_{\cdot}$ operators for the case of $|i-j|=1$.

\begin{definition}\label{remalpha-i}{\rm
% {\rm{(1)}}
 Let $e_i(\boldalpha)\simeq 0$ if $e_i$ attempts to act on a simple root $\alpha_i$ in the Kostant partition $\boldalpha$, i.e., if $e_i(\boldalpha)\ne\widehat{e}_i(\boldalpha)$, where necessarily $e_i(\boldalpha)=0$. Equivalently, let  $e_i(\boldalpha)\not\simeq 0$ if $e_i(\boldalpha)=\widehat{e}_i(\boldalpha)$, where the latter can be 0.
%{\rm{(2)}}
 Similarly, we let $e_ie_j\ldots e_l(\boldalpha)\simeq 0$ if $\bm{\beta}:=e_j\ldots e_l(\boldalpha)\ne 0$ and $e_i({\bm{\beta}})\simeq 0$.}
\end{definition}

%\vspace{12pt} The following are exceptions to the Stembridge axioms that are found when applying $e_{\cdot}$ operators for the case of $|i-j|=1$.

\vspace{6pt}
{\bf Exception BD} (broken diamond). We have 
$e_ie_j(\boldalpha)\simeq 0$ and  $e_je_i(\boldalpha)\simeq 0\,.$
 This happens only if $\widehat{\varepsilon}_i(e_j(\boldalpha))=\widehat{\varepsilon}_i(\boldalpha)$ and $ \widehat{\varepsilon}_j(e_i(\boldalpha))=\widehat{\varepsilon}_j(\boldalpha)\,.$

\vspace{12pt}{\bf Exception BO1} (broken octagon of first type). We have
$e_ie_j(\boldalpha)\simeq 0 $ and $ e_ie_j^2e_i(\boldalpha)\simeq 0\,.$
 This happens only if $\widehat{\varepsilon}_j(e_i(\boldalpha))=\widehat{\varepsilon}_j(\boldalpha)+1.$

\vspace{12pt}{\bf Exception BO2} (broken octagon of second type). We have
$e_i^2e_j(\boldalpha)\simeq 0 $ and $ e_ie_j^2e_i(\boldalpha)\simeq 0\,.$
This happens only if 
$\widehat{\varepsilon}_i(e_j(\boldalpha))=\widehat{\varepsilon}_i(\boldalpha)+1 $ and $\widehat{\varepsilon}_j(e_i(\boldalpha))=\widehat{\varepsilon}_j(\boldalpha)+1\,.$

\begin{figure}[H]
\centering
\includegraphics[scale=.35]{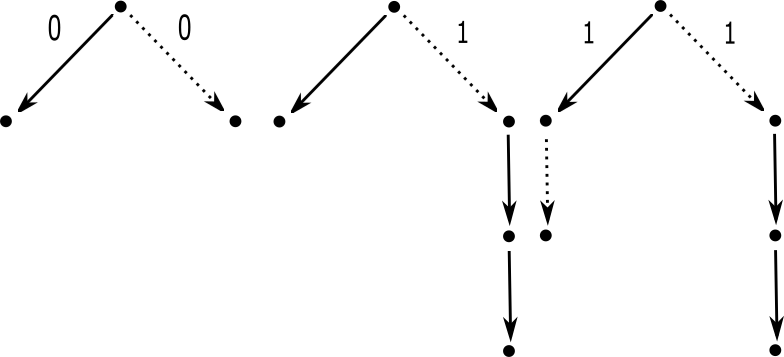}
\caption{Visualizations of the Stembridge Exceptions BD, BO1 and BO2 (left to right).  Here, the dotted arrows represent $e_i$ and solid arrows represent $e_j$.  The labels $0$ and $1$ denote respective changes in $\varepsilon_i$ on solid lines and $\varepsilon_j$ on dotted lines. } \label{Stembridge_exceptions_figure}
\end{figure}

%old theorems where labeled ije1

\begin{theorem}\label{stembridge ij>1} %\label{modified operators follow stembridge}
The Stembridge Axioms {\rm{S1}}, {\rm{S2}},  {\rm{S1}$'$}, and {\rm{S2}$'$} hold for the $e_\cdot$ and $f_\cdot$ operators, respectively, when $|i-j|>1$, where in the first two axioms we can use both the $\varepsilon_\cdot$ and $\widehat{\varepsilon}_\cdot$ functions.
\end{theorem}

\begin{theorem}\label{stembridge ij=1}  The Stembridge Axioms {\rm{S1}$'$}$-${\rm{S3}$'$} hold for the $f_\cdot$ operators when $|i-j|=1$.
 Consider $i,j$ with $|i-j|=1$ and $\boldalpha$ with $e_i(\boldalpha)\ne 0$. The Stembridge Axiom {\rm S1} holds with $\varepsilon_j$ replaced with $\widehat{\varepsilon}_{j}$. If, in addition, $e_j(\boldalpha)\ne 0$, then we are in one of the situations described in the Stembridge Axioms {\rm S2} and {\rm S3} (we again use $\widehat{\varepsilon}_{\cdot}$), unless we have Exceptions {\rm BD}, {\rm BO1}, or {\rm BO2}, possibly with $i,j$ switched in the last two.
 
\end{theorem}

  \section{The $\phi_i$ Operator on Kostant Partitions} 
  
 We now use our $f_i$ and $e_i$ operators, which act on Kostant partitions and preserve their size, to define operators which act on pairs $(w,\boldalpha)\in S_{\lambda,\mu}$.  The following operators are, by their very nature, sign-reversing with respect to $w$ and size preserving with respect to $\boldalpha$.  We will then be able to use these operators to define an involution on the elements of $S_{\lambda,\mu}$ which will give the desired cancelation in the related Kostka-Foulkes polynomial.
  
\begin{definition}\label{phi def}
{\rm
For $i\in I$ and $k_i = \langle w(\lambda + \rho ),\alpha_i^{\vee}\rangle$, $\boldalpha\in\mathcal{K}$, define the operator $$\phi_i(\boldalpha) := \begin{cases} 
      f_i^{-k_i}(\boldalpha) & \text{if}\hspace{5pt} k_i<0 \\
      e_i^{k_i}(\boldalpha) & \text{if}\hspace{5pt} k_i>0 \hspace{5pt}\text{and}\hspace{5pt} e_i^{k_i}(\boldalpha)\neq 0 \\
      0 & \text{otherwise.}
   \end{cases}$$ Note that $k_i\neq 0$, as $\lambda + \rho$ is in the dominant cone and so $w(\lambda + \rho)$ does not lie on any reflecting hyperplane. 
   %\color{red}Further, we have $\phi_i(\boldalpha)\neq 0$ when $k_i<0$ (this will be proven in Lemma~\ref{Lecouvey f exists})\color{blue} do we need this? \color{black}
   Extend this operator to pairs $(w,\boldalpha)$ where $w\in S_n$ and $\boldalpha$ is a Kostant partition with weight $w(\lambda + \rho) - (\mu + \rho)$. This then gives way to a graph where there is an edge $(w, \boldalpha)\xrightarrow{i}(w\textprime{},\boldsymbol\beta)$ if and only if $w\textprime{} = s_iw$ and $\boldsymbol\beta = \phi_i(\boldalpha)= f_i^{-k_i}(\boldalpha)$.
   }
\end{definition}

Before we can relate these new operators to the Kostka-Foulkes polynomials, we need to understand the graph structure induced by the $\phi_i$ operators.  We now show that the $\phi_i$ operators satisfy the braid relations, and further, that the graphs induced by the $\phi_i$ operators actually consists of certain subgraphs of Cayley graphs. % The proof for Theorems~\ref{stembridge ij=1} and~\ref{stembridge ij>1} (pertaining to the Stembridge axioms) is given in Section~\ref{section stembridge proof} and the proofs for Theorems~\ref{braid relations f} and~\ref{braid relations e} (pertaining to the braid relations) are given in Section~\ref{section braid relation}.  
%The proof for Theorem~\ref{psi is involution} is given in Section~\ref{matching proof section}.
%We first provide some exceptions to the Stembridge axioms.

% The following shows that the $\phi_i$ operators satisfy the braid relations.

\begin{theorem}\label{braid relations f}
Suppose that for a given Kostant partition $\boldalpha$, we have $i,j$ such that $k_i,k_j<0$ and $\phi_i(\boldalpha)\neq 0$ and $\phi_j(\boldalpha)\neq 0$.  Then 

\begin{enumerate}
\item if $|i-j|>1$, $\phi_i\phi_j(\boldalpha) = \phi_j\phi_i(\boldalpha)\neq 0$ and
\item if $|i-j|=1$, $\phi_j\phi_i\phi_j(\boldalpha) = \phi_i\phi_j\phi_i(\boldalpha)\neq 0$.
\end{enumerate} 
\end{theorem}

%\vspace{12pt} The following theorem is conditional, based on the validity of Conjecture~\ref{no broken diamond in hexagon}.

\begin{theorem}\label{braid relations e}
Suppose that for a given Kostant partition $\boldalpha$, we have $i,j$ such that $k_i,k_j>0$ and $\phi_i(\boldalpha)\neq 0$ and $\phi_j(\boldalpha)\neq 0$.

\begin{enumerate}
\item If $|i-j|>1$, then $\phi_i\phi_j(\boldalpha) = \phi_j\phi_i(\boldalpha)\neq 0$.
\item Suppose that $|i-j|=1$.
\begin{enumerate}
\item If $\varepsilon_i(e_j^{k_j}\boldalpha)< k_i+k_j$ , then $\phi_j\phi_i\phi_j(\boldalpha) = \phi_i\phi_j\phi_i(\boldalpha) =\phi_i\phi_j(\boldalpha)=0$ but $\phi_j\phi_i(\boldalpha)\neq 0$.
\item Otherwise,  $\phi_j\phi_i\phi_j(\boldalpha) = \phi_i\phi_j\phi_i(\boldalpha)\neq 0$.
\end{enumerate}
\end{enumerate} 
\end{theorem}

\begin{remark}\label{Remark broken by}{\rm We will refer to the case of $2a$ in Theorem~\ref{braid relations e} as a \textit{broken hexagon} and we will say $\boldalpha$ \textit{is broken by} $\phi_i$ \textit{after} $\phi_j$.}
\end{remark}
The proof of Theorem~\ref{braid relations f} is a constructive one, using the Stembridge axioms to build larger, and more general, relations.  The proof of Theorem~\ref{braid relations e} is contingent upon the following conjecture which is believed to be true through computer testing.

\begin{conjecture}\label{no broken diamond in hexagon}
Let $\boldalpha\in S_{\lambda,\mu}$ and $k_i,k_{i+1}>0$ be defined as in Definition~\ref{phi def}.  Suppose that $\phi_i(\boldalpha)$ and $\phi_{i+1}(\boldalpha)$ are both non-zero.  Then the paths $E_i:=e_i^{k_{i+1}}e_{i+1}^{k_i+k_{i+1}}e_i^{k_i}(\boldalpha)$ and $E_{i+1}:=e_{i+1}^{k_{i}}e_{i}^{k_i+k_{i+1}}e_{i+1}^{k_{i+1}}(\boldalpha)$ in the graph induced by the modified crystal operators do not contain any edges that give way to Stembridge Exception \rm{BD}.
\end{conjecture}

 %\begin{remark}\label{bicrystal and wt remark}
 
%We note the following:
%\begin{enumerate}
% \item The resulting graph is a subgraph of the bi-crystal given by $B(\infty)$ and $B^*(\infty)$ where the  $f_i$ operators can be explicitly stated in terms of $\widetilde{f}_i$ and $\widetilde{f}^*_i$ via the following:

 %  $$f_i(\alpha) = \begin{cases}
 %    \widetilde{f}_i(\alpha) & \text{if}\hspace{5pt} b_i^* <  a_i  \\
 %    \widetilde{f}_i^*(\alpha) & \text{if} \hspace{5pt} b_i^*\geq  a_i \hspace{3pt} \text{and} \hspace{5pt} a^*_i>0  \\
  %   0                     & \text{otherwise}
  %\end{cases}.$$
  
  %While this is a nice observation, we will refrain from depending on this structure since corresponding combinatorial models are not as readily available in other Lie types.
  
%\item If we let $wt(\alpha)$ be the weight of $\alpha$ when regarded as a vertex of $M(\infty)$ the set of all multisegments (or Kostant partitions), then $wt(\alpha) = -\overline{\alpha}$.  Further, if we let $wt_i(\alpha) = \langle wt(\alpha),\alpha_i^{\vee}\rangle$, then we have $wt_i(\alpha) = a_i-b_i+a^*_i-b^*_i = u_i - v_i$.
  
 % \end{enumerate}
 %\end{remark} 
 
% \section{The Phi Operator and Braid Relations}

% \begin{theorem}\label{psi is involution}
% The function $\psi$ is a sign changing involution on $S_{\lambda,\mu}$
% \end{theorem}

%We can now prove the main theorem.

%\subsection{An Involution on $S_{\lambda,\mu}$ and the Kostka-Foulkes Polynomials}
\section{An Involution on $S_{\lambda,\mu}$ and the Main Result}\label{Section Main Result}

   The graphs resulting from the $\phi_i$ operators acting on $S_{\lambda,\mu}$ have nice additional structures which we will now discuss.  This structure will then give way to the main theorem.
   
   \begin{theorem}\label{subgraph of cayley}
Each connected component of the graph given by the $\phi_i$ operators is a subgraph $\Gamma$ of the Cayley graph of $S_n$ such that $w\in\Gamma\implies u\in\Gamma$ for any $u\leq w$ in the weak Bruhat order.
\end{theorem}

We note that for $(id,\boldalpha)$ in $S_{\lambda,\mu}$ we have $k_i = \langle\lambda + \rho,\boldalpha_i^{\vee}\rangle>0$  for $i\in\{1,\hdots,n-1\}$.

\begin{proposition}\label{all edges on top}
Assume $i\in\{1,\hdots,n-1\}$ is such that $e_i^{k_i}(\boldalpha)=0.$  Then for any $(w,\boldalpha')$ in the connected component of $(id,\boldalpha)$ where $k'_i>0$, we have $e_i^{k'_i}(\boldalpha) = 0$.  Equivalently, in the component of $(id,\boldalpha)$, an arrow $\xrightarrow{i}$ can exist only if there is such an arrow starting from $(id,\boldalpha)$.
\end{proposition}

%\begin{proof}
%Assume that there exists $(w,\boldalpha')$ in the connected component $C$ of $(id,\boldalpha)$ such that $e_i^{k'_i}(\boldalpha)\neq 0$.  One can assume that $l(w)$ is minimal along all those possible $(w,\boldalpha')$.  In particular, $l(w)\geq 1$.  Since $l(w)$ is minimal, $w$ belongs to the parabolic subgroup of $S_n$ generated by the $s_j$ with $i\neq j$.  Otherwise, $s_i$ appears in any minimal decomposition of $w$. Since $(w,\boldalpha')$ does belong to $C$, there should exist at least one such decomposition of $w=us_iv$ yielding an arrow $(v,\beta)\xrightarrow{i} (s_iv,\beta')$ in $C$ with $l(v)<l(w)$, thus a contradiction.  We so derive that $\boldalpha'$ is obtained from $\boldalpha$ by applying operators $e_j$ with $j\neq i$.  Let us write $w_i^r = (+)^{u_i}(-)^{v_i}$ and $(w')_i^r = (+)^{u'_i}(-)^{v'_i}$ for $\boldalpha$ and $\boldalpha'$.  On the one hand, by Remark~\ref{G_i graph change}, we have that $u'_i\leq u_i$.  On the other hand, we have $k_i>0$ and $k'_i = \langle w(\lambda + \rho),\boldalpha_i^{\vee}\rangle>0$ since $w$ belongs to the parabolic subgroup generateLemmad by the $s_j$ with $j\neq i$.  We can thus apply assertion $3$ of Lemma~\ref{Lecouvey f exists} and get $u_i < \langle\mu + \rho,\boldalpha_i^{\vee}\rangle\leq u'_i$.  This gives the desired contradiction.
%\end{proof}

\begin{corollary}\label{subgraph of cayley of parabolic subgroup}
Each connected component $C$ of the graph given by the $\phi_i$ operators is isomorphic to a subgraph $\Gamma$ of the Cayley graph of a parabolic subgroup of $S_n$ determined by its source vertex such that $w\in\Gamma\implies u\in\Gamma$ for any $u\leq w$.
\end{corollary}

%\begin{proof}
%This follows from Theorem~\ref{braid relations f}, Theorem~\ref{subgraph of cayley} and Proposition~\ref{all edges on top}.
%\end{proof}

We can conjecturally describe these subgraphs further via the following which gives a larger scale version of Proposition $\ref{all edges on top}$.

\begin{conjecture}\label{broken hexagon propogates up}
If in a given component $C$ of the graph given by the $\phi_i$ operators there is a broken hexagon broken by some $\phi_i$ after some $\phi_j$, then the source vertex of the component also has a broken hexagon broken by $\phi_i$ after some $\phi_j$. This has been confirmed via computer testing.% on Sage.
\end{conjecture}

%\begin{remark}
%This has been confirmed via computer testing on Sage.
%\end{remark}

We now define a matching on each component $C$ of the graph given by the $\phi_i$ operators.  We will show that the following function is actually an involution.

\begin{definition}\label{psi def} {\rm We define a function $\psi:S_{\lambda,\mu}\rightarrow S_{\lambda,\mu}$.  Let $C$ be a connected component of the graph induced on $S_{\lambda,\mu}$ by the $\phi_i$ operators with source vertex $(id,\boldalpha)$. Then $C \cong S_{n_1}\times S_{n_2}\times\hdots\times S_{n_r}$, where $n_1\leq n_2\leq \hdots \leq n_r$ and each $n_k$ corresponds to the length of the $k$-th block of sucessive indices of the arrows in $C$ (see Corollary~\ref{subgraph of cayley of parabolic subgroup}). Let $\psi_C$ be the restriction of $\psi$ to the component $C$.  If $C$ consists only of $\boldalpha$, then let $\psi_C$ be the identity map.  Otherwise, let $\psi_C = \phi_i$ where $i$ is determined as follows:

\begin{enumerate}
\item If $\boldalpha$ is not broken by $\phi_1$ after $\phi_2$, let $i=1$.
\item Otherwise, let $i = max\{j\in [n_1-1]:\boldalpha \hspace{4pt}\text{is broken by } \phi_{l-1} \text{ after } \phi_l \text{ for all } 1< l\leq j\}$ (cf. Remark~\ref{Remark broken by}).
\end{enumerate}
}
 \end{definition}
 
Note that Definition~\ref{psi def} is based on viewing labels on $C$ through the isomorphism with subgraphs of Cayley graphs of symmetric groups.  This means that the labels on the actual graph $C$ may be translated and not actually start at $1$.%  See Example~\ref{inv examples}~(2) where the graphs are isomorphic to subgraphs of $S_3$ but their edge labels are $2$ and $3$ instead of $1$ and $2$.
 
%\vspace{12pt}
 We can now state our main theorem.  Note that, for the remainder of this section, \textbf{we are assuming the validity of Conjectures~\ref{broken hexagon propogates up} and~\ref{no broken diamond in hexagon}}. 

\begin{theorem}\label{KF main theorem} We write the Kostka-Foulkes polynomial with all positive coefficients as follows:
$$K_{\lambda,\mu}(t) = \sum\limits_{(w,\boldalpha)\in S_{\lambda,\mu}}sgn(w)t^{|\boldalpha |} = \sum\limits_{(w,\boldalpha)\in D}sgn(w)t^{|\boldalpha |}= \sum\limits_{(id,\boldalpha)\in D}t^{|\boldalpha |} $$ where $D$ is the collection of fixed points of $\psi$. Furthermore, the terms of these fixed points correspond directly to the elements $\boldalpha\in S_{\lambda,\mu}$ where $\phi_i(\boldalpha)=0$ for all $i\in I$.
\end{theorem}

\begin{figure}[H]
\centering
\includegraphics[scale=.45]{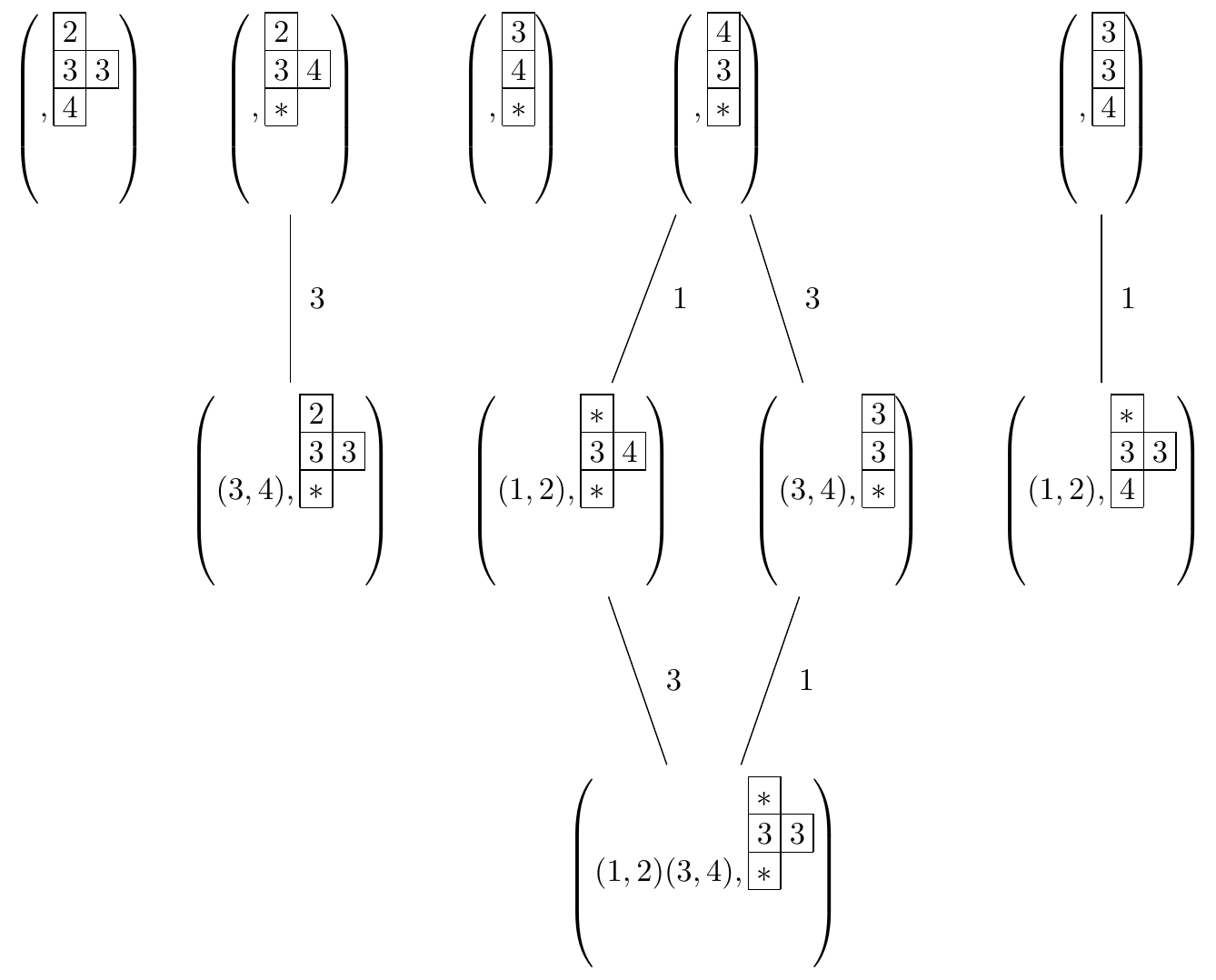}
\caption{The graph induced by the $\phi_{\cdot}$ operators acting on vertices from $S_{[2,2],[1,1,1,1]}$.} \label{phi_graph_ex}
\end{figure}

\begin{example}\label{inv examples} 
{\rm
 We now consider an example using the tableau notation for Kostant partitions, $\boldalpha$, as was used in \cite{Lee Salisbury} and are to be read as follows: each entry $j$ in row $i$ represents the positive root $\epsilon_i - \epsilon_j$ in the formal sum that is $\boldalpha$. 

%\begin{enumerate}

%\item
 Refer to Figure~\ref{phi_graph_ex} for the graph induced by the $\phi_i$ operators acting on elements from $S_{[2,2],[1,1,1,1]}$ (cf. Example~\ref{example KF poly}).  We can also see the pairing given by $\psi$ on each component.  Note that the two singleton components (the fixed points of $\psi$) are of sizes $2$ and $4$, which give way to the reduced Kostka-Foulkes polynomial $t^4+t^2$ from Example~\ref{example KF poly}. Further, note that this example shows that we need the actions from both $B^*(\infty)$ and $B(\infty)$.

%\item We now look at an example showing the necessity of defining $\psi$ on individual connected components. Refer to Figure~\ref{diff_breaks} where we have two connected components, $C_1$ and $C_2$, of the graph induced by the $\phi_i$ operators acting on elements from $S_{[3,1,1,1,0,0,0],[1,1,1,1,1,1,0]}$. Note that, although both are isomorphic to subgraphs of $S_3$, in $(a)$ we have $\psi_{C_1} = \phi_3$ and in $(b)$ we have $\psi_{C_2} = \phi_2$.  

%\item 
%As stated in the introduction, the classical embedding of SSYT into $B(\infty)$ does not give way to the Kostka-Foulkes polynomials.  Consider  $\lambda=[2,2,0,0], \mu=[1,1,1,1]$, and $n=4$.  Then our current model gives the correct $K_{\lambda,\mu}(t)=t^2+t^4$, while the image of the classic crystal embedding gives $t^2+t^3$.

%\end{enumerate}

}
\end{example}

We now wish to describe the fixed points of $\psi$ in a more direct combinatorial manner.  Recall that the fixed points of $\psi$ were exactly the elements $\boldalpha\in S_{\lambda,\mu}$ where $\phi_i\boldalpha=0$ for all $i\in I$.  We also found that these fixed points were only of the form $(id, \boldalpha)$ and so $k_i =\langle \lambda+\rho,\alpha_i^{\vee}\rangle >0$ for all $i\in I$.     

%\begin{definition}
%Consider a Kostant partition $\boldalpha$ and the associated words $W_i^*W_i$ and positive integers $k_i$ for $i\in I$.  If, for any $i$, when looking at the word $W_i^*W_i$ from the right, there are never $k_i$ or more $-$'s than $+$'s, and, if there are, at least one of the $k_i$ left-most $-$'s are associated to a simple root $(i,i+1)$, we say that $\boldalpha$ is \textit{admissible}. 

%\end{definition}
\begin{definition} {\rm
Consider a Kostant partition $\boldalpha$ and the associated words $W_i^*W_i$ and positive integers $k_i$ for $i\in I$.   We say that $\boldalpha$ is \textit{admissible} if for any $i$, when looking at the word $W_i^*W_i$ from the right, there are
\begin{enumerate}
\item never $k_i$ or more $-$'s than $+$'s, or
\item at some point there are $k_i$ or more $-$'s than $+$'s but at least one of the $k_i$ left-most $-$'s are associated to a simple root $(i,i+1)$.
\end{enumerate}
}

\end{definition}

\begin{corollary}

Consider the set $S_{\lambda,\mu}^{w=id} := \{(id,\boldsymbol\alpha): \boldsymbol\alpha\in\mathcal{K}\hspace{4pt} \text{where} \hspace{4pt} \overline{\boldsymbol\alpha} = \lambda - \mu\}$. Then the collection of fixed points of $\psi$ is $D =\{ \boldalpha\in S_{\lambda,\mu}^{w=id}: \boldalpha \hspace{4pt}\text{is}\hspace{4pt} \text{admissible}\} $ and we again have
$$K_{\lambda,\mu}(t) = \sum\limits_{\boldalpha\in D}t^{|\boldalpha |}. $$

\end{corollary}

%\end{document}

%Here is a citation~\cite{greenwade93} with URL. Here is another
%citation~\cite{MR4}, which is the earliest-numbered Math Review with a
%functioning DOI: MR0000004. Finally, a reference to
%equation~\eqref{eqn:eq1}. And a reference to Figure~\ref{fig:plot}.

%\acknowledgements{\lipsum[20]}

%% if you use biblatex then this generates the bibliography
%% if you use some other method then remove this and do it your own way
%\printbibliography

\end{document}